\documentclass[12pt]{article}
  \usepackage[cp1251]{inputenc}
  \usepackage[english]{babel}
  \usepackage{amsfonts,amssymb,amsmath,amstext,amsbsy,amsopn,amscd,graphicx}
  \usepackage{graphics}
  \usepackage{amsthm}

  \def\thtext#1{
  \catcode`@=11
  \gdef\@thmcountersep{. #1}
  \catcode`@=12}

  \newtheorem{theorem}{Theorem}[section]
  \newtheorem{ass}{Statement}[section]
  \newtheorem{prop}{Proposition}[section]
  \newtheorem{lem}{Lemma}[section]
  \newtheorem{cor}{Corollary}[section]

 \newcounter{il}[section]
 \renewcommand{\thtext}
 {\thechapter.\arabic{il}}

 \newcounter{il2}[section]
 \renewcommand{\thtext}
 {\thechapter.\arabic{il2}}

 \newenvironment{dfn}{\trivlist \item[\hskip\labelsep{\bf Definition}]
 \refstepcounter{il}{\bf \arabic{section}.\arabic{il}.}}%
 {\endtrivlist}

 \newenvironment{rk}{\trivlist \item[\hskip\labelsep{\bf Remark}]
 \refstepcounter{il2}{\bf \arabic{section}.\arabic{il2}.}}%
 {\endtrivlist}

 \newenvironment{example}{\trivlist \item[\hskip\labelsep{\bf Example.}]}%
 {\endtrivlist}

 \catcode`@=11
 \def\.{.\spacefactor\@m}
\def\relaxnext@{\let\next\relax}
\def\nolimits@{\relaxnext@
 \DN@{\ifx\next\limits\DN@\limits{\nolimits}\else
  \let\next@\nolimits\fi\next@}%
 \FN@\next@}
\def\newmcodes@{\mathcode`\'39\mathcode`\*42\mathcode`\."613A%
 \mathcode`\-45\mathcode`\/47\mathcode`\:"603A\relax}
\def\operatorname#1{\mathop{\newmcodes@\kern\z@
 \operator@font#1}\nolimits@}
 \catcode`@=12
 \def\rom#1{{\em#1}}

 \def\rank{\operatorname{rank}}
 \def\corank{\operatorname{corank}}

 \def\Loc{\operatorname{Loc}}
 \def\Sym{\operatorname{Sym}}

 \def\0{{\mathbf 0}}
 \def\1{{\mathbf 1}}
 \def\a{{\mathbf a}}
 \def\b{{\mathbf b}}
 \def\c{{\mathbf c}}
 \def\u{{\mathbf u}}
 \def\v{{\mathbf v}}
 \def\w{{\mathbf w}}

\title {Framed 4-Valent Graphs: Euler Tours, Gauss Circuits and Rotating Circuits}

\author{D.\,P.~Ilyutko\footnote{Partially supported by grants of RF President
NSh -- 660.2008.1, RFFI 07--01--00648, RNP 2.1.1.3704, the Federal
Agency for Education NK-421P/108.}}
\date {}

\begin {document}

 \maketitle
 \abstract{In the present paper we give an explicit formula which
 allows us immediately to describe a unique Gauss circuit on a framed $4$-valent
 graph (a graph with a structure of opposite edges) from an arbitrary
 Euler tour on the graph whenever the Gauss circuit exists. This formula only depends on the
 adjacency matrix of an Euler tour and also tells us
 whether there exists a Gauss tour on a framed 4-valent graph or not.
 It turns out that the results are also valid for all symmetric matrices
 (not just realisable by a chord diagram).}


 \section {Introduction}


In the present paper, we consider finite connected {\em framed
4-graphs}, i.e.\ 4-valent graphs with a pairing of the (half)edges
emanating from each vertex. The two (half)edges belonging to the
same pair are said to be {\em opposite} at this vertex. An example
of a framed 4-graph is a diagram of a virtual link where classical
crossings play the role of vertices, and virtual crossings are just
intersection points of images of different edges, and edges forming
overpass or underpass are opposite to each other. Graphs obtained
from links in such a way are called {\em projections of links}.

Each framed 4-graph $G$ has an {\em Euler tour} $U$, i.e.\ a
continuous map from the circle $S^{1}$ onto $U$ and this map is
bijective outside the vertices of $G$ and has exactly two inverse
images at each vertex of $G$. In each vertex we have two
possibilities of running edges: we move from an (half)edge to the
opposite to it (half)edge; we move from an (half)edge to a
non-opposite to it (half)edge. There are two special sorts of Euler
tours on a given framed 4-graph $G$: only the first possibility
occurs in each vertex of an Euler tour, and only the second
possibility occurs in each vertex of an Euler tour. It is not
difficult to see that Euler tours of the first sort do not always
exist and if they do then there exists only a unique one (the {\em
Gauss circuit}), but Euler tours of the second sort ({\em rotating
circuits}) exist on any connected framed 4-graphs and the number of
them for a 4-graph is larger than one, see~\cite {IM1,IM2,Ma1}. If
we consider a projection of a knot, i.e. a link with one component,
then it has the Gauss circuit, but a projection of a link with
larger than $1$ components does not have a Gauss circuit.

In low-dimensional topology both approaches, the Gauss circuit
approach and the rotating circuit approach, are very widely used.
The Gauss circuit approach is applied in knot theory, namely in the
construction of finite-type invariants, Vassiliev invariants~\cite
{BN,GPV,CDL}, and in the planarity problem of immersed curves,~\cite
{CE1,CE2,RR}. However, for detecting planarity of a framed 4-graph
it is more convenient to use the rotating circuit approach,
see~\cite {Ma1,Ma2}. The criterion of the planarity of an immersed
curve, which is the framed 4-graph, is formulated very easy: an
immersed curve is planar if and only if the chord diagram obtained
from a rotating circuit is a {\em d-diagram}, i.e.\ the set of all
the chords can be split into two sets and the chords from one set do
not intersect each other, see~\cite {Ma3}. If we want to generalise
the planarity problem to the problem of finding the minimum genus of
a closed surface which a given curve can be immersed in the rotating
circuit approach is also more useful. There are criteria giving us
the answer to the question what is the minimum genus for a given
curve, see~\cite {Ma1}.

Since there are many rotating circuits corresponding to the same
Gauss circuit many properties of the Gauss circuit can be read out
of any of these rotating circuits no matter which one is considered.
Consequently, these properties do not depend on the particular
choice of a rotating circuit. For instance, if one of rotating
circuits has a d-diagram then the other ones do the same. Thus it is
necessary to obtain an easy formula allowing us to get the Gauss
circuit from a rotating circuit and vice versa. Of course having a
framed 4-graph we can answer the question whether there is or not a
Gauss circuit on it and find the Gauss circuit whenever it exists by
just traveling along our 4-graph. But this method does not reflect
explicit relations between topology and combinatorics of Euler tours
if we have a 4-graph with many vertices. In the present paper we
give an explicit formula which depends only on the adjacency matrix
of an Euler tour: taking any Euler tour and constructing its
adjacency matrix we can find the adjacency matrix of the Gauss
circuit. It turns out that the given formula is also valid for all
symmetric matrices (not just realisable by a chord diagram).
Investigating this formula we can get some interesting facts of
symmetric matrices.

The present paper is organised as follows.

We first give main definitions concerning 4-graphs, framed 4-graphs,
Euler tours, and preliminary results about them.

In the second section we investigate the question of the existence
of a Gauss circuit in terms of its adjacency matrix, and in the
third section an explicit formula connecting adjacency matrices of
distinct Euler tours is given. Using this formula we can easy get
the adjacency matrix of the Gauss circuit.

In the fourth section we generalise the results of the preceding
sections to the case of all symmetric matrices not only matrices
realised by chord diagrams.


 \section {Main Definitions and Preliminary Results}


During the whole article by graph we mean a connected finite graph,
possibly, having loops and/or multiple edges.

\subsection {4-Valent Graphs and Euler Tours}

Let $G$ be a graph with the set of vertices $V(G)$ and the set of
edges $E(G)$. We think of an edge as an equivalence class of the two
half-edges. We say that a vertex $v\in V(G)$ has the {\em degree}
$k$ if $v$ is incident to $k$ half-edges. A graph whose vertices
have the same degree $k$ is called {\em $k$-valent} or a {\em
k-graph}. The free loop, i.e. the graph without vertices, is also
considered as $k$-valent graph for any $k$.

Let $H$ be a connected 4-graph on the set of vertices $V(H)$ and let
$U$ be an Euler tour of $H$, i.e. a tour while traveling along it we
run each edge exactly one time. For every vertex $v\in V(H)$ there
are precisely two closed paths $P_v$ and $Q_v$ on $U$ having no
common edges, starting and ending at $v$ and having no internal
vertex equal to $v$. There exists precisely one Euler tour distinct
from $U$ also connecting the paths $P_v$ and $Q_v$ (if we run along
$U$ in some direction then in the new Euler tour we run along $P_v$
according to the orientation of $U$, and run along $Q_v$ according
to the reverse orientation of $U$). Let us denote by $U*v$ the new
Euler tour obtained from $U$. The transformation $U\mapsto U*v$ has
been introduced by Kotzig~\cite {Ko} who called it a {\em
$k$-transformation}. He proved the following statement.

 \begin {prop}
Any two Euler tours of a $4$-graph are related by a sequence of
$k$-transformations.
 \end {prop}

Let $w=x_1x_2\ldots x_{k-1}x_k$ be a word, i.e.\ a sequence of
letters from some finite alphabet $X$. The {\em mirror image of $w$}
is $\widetilde{w}=x_kx_{k-1}\ldots x_2x_1$. We will consider the
class of words where each word from this class is either a cyclic
permutation $w_i=x_ix_{i+1}\ldots x_kx_1\ldots x_{i-1}$, $1\leqslant
i\leqslant k$, of $x_1\ldots x_k$ or the mirror image of a cyclic
permutation $w_i$. We denote this class by $(x_1\ldots x_k)$ and we
call this class a {\em cyclic word}.

 \begin {dfn}
A word is called a {\em double occurrence word} if each its letter
occurs twice in it.
 \end {dfn}

Let $X$ be a finite set. Let $m$ be a double occurrence cyclic word
over $X$, i.e. a class of words. Then $m$ has a chord
representation, which is constructed by placing successively the
letters of $m$ around a circle $S^1$, choosing a point of $S^1$ near
each occurrence of a letter and joining by a chord each pair of
points corresponding to the two occurrences of the same letter. It
is not difficult to see that we get the one-to-one correspondence
between the set of double occurrence cyclic words and the set of
chord diagrams.

\begin {figure} \centering\includegraphics[width=150pt]{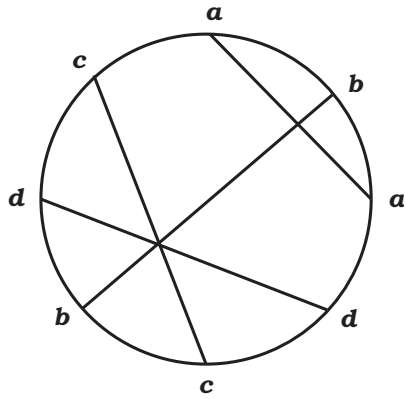}
  \caption{A Chord Representation of $(abacdbcd)$} \label{chre}
 \end {figure}

 \begin {example}
Consider $m=(abacdbcd)$. The word $m$ has the chord representation
depicted in Fig.~\ref {chre}.
 \end {example}

\begin {figure} \centering\includegraphics[width=250pt]{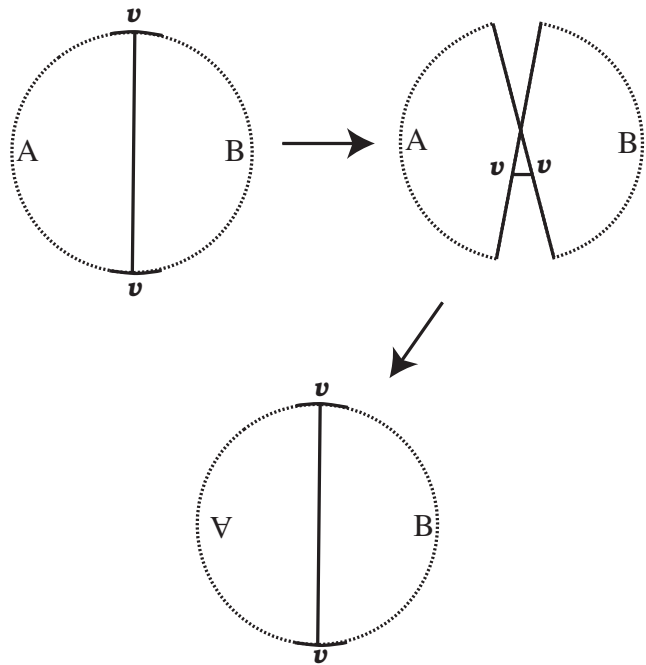}
  \caption{The Operation $*$ on Chord Diagrams} \label{staronch}
 \end {figure}

Define the operation $*$ on words which will correspond to the
$k$-transformation. Let $m=(vAvB)$ where $A,\,B$ are subwords of
$m$, and letters belong to some finite alphabet. Then we define
$m*v=(v\widetilde{A}vB)$, $\widetilde{A}$ is the mirror image of
$A$. In Fig.~\ref {staronch} the transformation $m\mapsto m*v$ is
depicted for chord diagrams (dashed arcs of chord diagrams contain
the ends of all the chords distinct from $v$). Mostly for each
transformation on a chord diagram we assume that only a fixed
fragment of the chord diagram is being operated on. The pieces of
the chord diagram not containing chords participating in this
transformation are depicted by dashed arcs.

Let $U$ be an Euler tour of a connected $4$-graph $H$ with the set
of vertices $V(H)=\{v_1,\dots,v_n\}$, which is also considered as an
alphabet. When traveling along $U$ we meet each vertex twice. Let us
denote by $m(U)$ the (cyclic) word over $V(H)$ that equals (up to
cyclic equivalence) the sequence of the vertices that are
successively met along $U$. It is obvious that in the obtained word
each vertex appears precisely twice then Euler tours are encoded by
double occurrence cyclic words. It follows from the definition that
$m(U*v)=m(U)*v$ and if we have a double occurrence cyclic word $m$
or a chord diagram we can construct the $4$-graph having an Euler
tour $U$ such that $m(U)=m$. We just contract each pair of vertices
of the chord diagram labeled by a same letter (a chord) into a
single vertex and identify the new vertex with this letter.

\subsection {Framed 4-Valent Graphs and Euler Tours}

In the present subsection we will consider $4$-graphs with the
additional structure.

 \begin {dfn}
A $4$-graph is called {\em framed} if for every vertex the four
emanating half-edges are split into two pairs of (formally) opposite
edges. The edges from one pair are called {\em opposite to each
other}.
 \end {dfn}

\begin {figure} \centering\includegraphics[width=250pt]{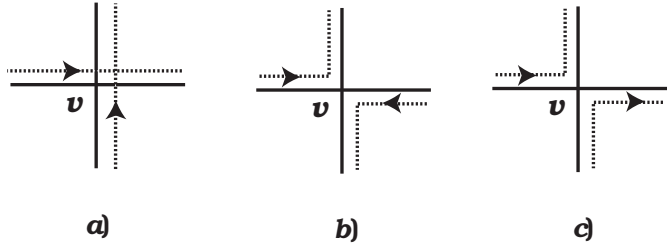}
  \caption{Passing through the Vertex} \label{passthver}
 \end {figure}

Let $H$ be a framed $4$-graph and $U$ be an Euler tour on $H$.
Construct the framed double occurrence cyclic word $m(U)$ (the
framed chord diagram) corresponding to $U$. In each vertex $v$ of
$H$ we have the following three possibilities of running along $U$
through $v$:
 \begin {enumerate}
  \item
we pass from a half-edge to the opposite to it half-edge, see
Fig.~\ref {passthver}~a). The vertex $v$ is called a {\em Gaussian
vertex for $U$} and the chord corresponding to this vertex is also
called a {\em Gaussian chord};
  \item
we pass from a half-edge to a non-opposite to it half-edge, and the
orientations of opposite edges are different, see Fig.~\ref
{passthver}~b). The vertex $v$ is called a {\em non-Gaussian vertex
for $U$ with framing $0$} and the chord corresponding to this vertex
is also called a {\em non-Gaussian chord with framing $0$};
  \item
we pass from a half-edge to a non-opposite to it half-edge, and the
opposite edges have the same orientation, see Fig.~\ref
{passthver}~c). The vertex $v$ is called a {\em non-Gaussian vertex
for $U$ with framing $1$} and the chord corresponding to this vertex
is also called a {\em non-Gaussian chord with framing $1$}.
 \end {enumerate}

 \begin {dfn}
An Euler tour having only Gaussian vertices is called a {\em Gauss
circuit}. An Euler tour having only non-Gaussian vertices is called
a {\em rotating circuit} (see~\cite {IM1,IM2,Ma1,Ma2}).
 \end {dfn}

When running along an Euler tour $U$ on a 4-graph $H$ we meet each
vertex of $H$ twice. Now we are ready to construct the framed double
occurrence cyclic word $m(U)$ corresponding to $U$. Words will be
constructed over the alphabet $X=V(H)\cup V(H)^{-1}\cup V(H)^G$,
where $V(H)^{-1}$ is the set of letters $v^{-1}$ for each $v\in
V(H)$ and $V(H)^G$ is the set of letters $v^G$ for each $v\in V(H)$.
Vertices (of $H$) of the first type will be denoted in $m(U)$ by the
same symbols as in $V(H)$ but with the superscripts $G$, i.e.\ by
letters from $V(H)^G$. For example, $m(U)=(Av^G Bv^G)$ if $v$ is the
vertex of the first type. Vertices (of $H$) of the second type will
be denoted by the same symbols as in $V(H)$ but with the
superscripts $\{\pm1\}$, and the superscripts are the same for both
occurrences of the same letter, i.e.\ by letters from $V(H)\cup
V(H)^{-1}$. For example, $m(U)=(AvBv)$ (in practice, the
superscripts $+1$ are omitted) or $m(U)=(Av^{-1}Bv^{-1})$ if $v$ is
the vertex of the second type (we will not make a difference between
these words, i.e.\ we consider an equivalence class). Vertices of
the third type will be denoted by the same symbols as in $V(H)$ but
with the superscripts $\{\pm1\}$, and the superscripts are
different, i.e.\ by letters from $V(H)\cup V(H)^{-1}$. For example,
$m(U)=(AvBv^{-1})$ or $m(U)=(Av^{-1} Bv)$ if $v$ is the vertex of
the third type (we will not make a difference between these words,
i.e.\ we consider an equivalence class).

Depicting a double occurrence cyclic word by a chord diagram we will
use thick chords for vertices with framing $0$, dashed chords for
vertices with framing $1$, and chords with the label $G$ for
Gaussian vertices.

\begin {figure} \centering\includegraphics[width=150pt]{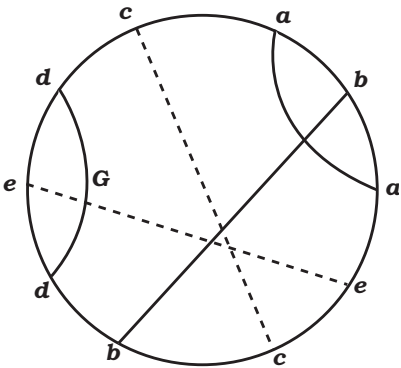}
  \caption{A Framed Chord Diagram for $(ab^{-1}acd^Ge^{-1}d^Gb^{-1}c^{-1}e)$} \label{frchdi}
 \end {figure}

 \begin {example}
Consider the framed word $m=(ab^{-1}acd^Ge^{-1}d^Gb^{-1}c^{-1}e)$.
We have: $d$ is a Gaussian letter, $a,\,b$ are non-Gaussian letters
with framing $0$ and $c,\,e$ are non-Gaussian letters with framing
$1$. The corresponding framed chord diagram is depicted in Fig.~\ref
{frchdi}.
 \end {example}

Let $V$ be a finite set. Having a framed double occurrence cyclic
word (a framed chord diagram) $m$ over $V\cup V^{-1}\cup V^G$ we can
construct the framed $4$-graph having an Euler tour $U$ such that
the framed word $m(U)$ coincides with $m$. We construct the
$4$-graph and then define the type of each vertex.

 \begin {rk}
When we consider framed double occurrence cyclic words over an
alphabet it is only important for us the positions of the same
letters but not its symbols, see~\cite {Tu}.
 \end {rk}

Let us define the {\em framed star} operation on the set of framed
double occurrence words. We denote this operation by the symbol $*$.

 \begin {rk}
We use the same notations as for double occurrence cyclic words.
Further we will consider only framed double occurrence cyclic words
and the same notations will not cause confusion.
 \end {rk}

\begin {figure} \centering\includegraphics[width=250pt]{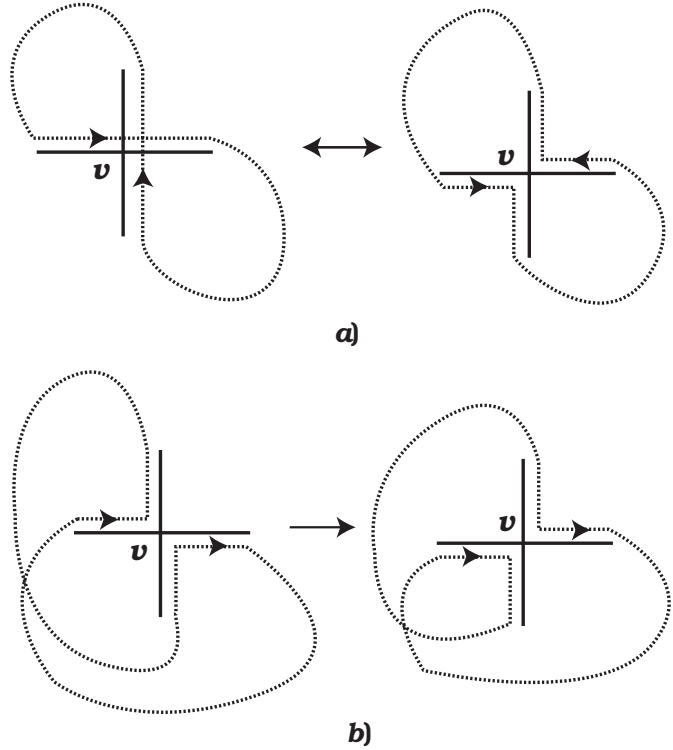}
  \caption{The Framed Star Operation} \label{croci}
 \end {figure}

Firstly, we construct the operation $\overline{w}$, where $m$ is an
arbitrary subword (not necessarily a double occurrence word) of a
framed double occurrence cyclic word. Let
$w=x^{\varepsilon_1}_1\dots x^{\varepsilon_k}_k$. Then
$\overline{w}=\overline{x^{\varepsilon_k}_k}\dots
\overline{x^{\varepsilon_1}_1}$, here
$\overline{x^{\varepsilon_l}_l}=x^{\varepsilon_l}_l$ if
$\varepsilon_l=g$,
$\overline{a^{\varepsilon_l}_l}=a^{-\varepsilon_l}_l$ if
$\varepsilon_l=\pm1$. Further, $m=(a^{\varepsilon}
m_1a^{\varepsilon'}m_2)$ is a double occurrence cyclic word. We have
$m*a=(a\overline{m_1}am_2)$ if $\varepsilon=\varepsilon'=g$
(Fig.~\ref {croci} a)); $m*a=(a^g\overline{m_1}a^gm_2)$ if
$\varepsilon=\varepsilon'\ne g$ (Fig.~\ref {croci} a));
$m*a=(a\overline{m_1}a^{-1}m_2)$ if $\varepsilon=-\varepsilon'$
(Fig.~\ref {croci} b)), i.e.\ applying the framed star to a Gaussian
letter it is transformed to the non-Gaussian letter with framing
$0$, applying to a non-Gaussian letter with framing $0$ it is
transformed to the Gaussian letter and applying to a non-Gaussian
letter with framing $1$ it is transformed to the non-Gaussian letter
with the same framing $1$.

 \begin {ass}[\cite{IM1,IM2,Ma1}]\label {ass:frwo}
Any two framed letters obtained from a framed $4$-graph are related
to each other by a sequence of the framed star operations.
 \end {ass}

The following corollary immediately follows from Statement~\ref
{ass:frwo}. We just apply the framed star operations to Gaussian
letters.

 \begin {cor}\label {cor:ro_ci}
Every framed $4$-graph has a rotating circuit.
 \end {cor}

 \begin {rk}
It is not difficult to prove Corollary~\ref {cor:ro_ci} by using
other methods, but we want to get used to the framed star operation.
 \end {rk}

It is obvious that there are many rotating circuits and that not
every framed $4$-graph has a Gauss circuit (if it has a Gauss
circuit then this tour is unique). The next section tells us whether
there exists a Gauss circuit or not, and how to get it whenever it
exists. The following theorem tells us how two rotating circuits are
related.

 \begin {ass}[\cite{IM1,IM2,Ma1}]\label {ass:2roci}
Any two rotating circuits given by framed double occurrence cyclic
words are related by a sequence of the following two
operations{\em:} the framed star operation applying to a
non-Gaussian letter with framing $1$, and $(((m*a)*b)*a)$, here $m$
is a framed double occurrence cyclic word, $a,\,b$ are non-Gaussian
letters with framing $0$ and they alternate in $m$, i.e.\ $m=(\ldots
a\ldots b\ldots a\ldots b)$.
 \end {ass}


 \section {Gauss Circuits}


In this section we give an explicit formula which allows us
immediately to describe a unique Gauss circuit on a framed $4$-graph
from an arbitrary Euler tour on the graph whenever the Gauss circuit
exists.

\subsection {The Existence of a Gauss Circuit}

We shall need two notions for establishing a criterion of the
existence of a Gauss circuit: the {\em adjacency matrix} of a framed
double occurrence cyclic word (a framed chord diagram) and a {\em
surgery} along chords.

 \begin {dfn}
Let $D$ be a framed chord diagram. We call two chords {\em linked}
if the ends of the first chord belong to the different connected
components of $S^1$ without the second chord's ends. Using the
language of framed double occurrence cyclic words we call two
letters $a,\,b$ {\em alternate} if we meet them alternatively
$(\ldots a^{\alpha_1}\ldots b^{\beta_1}\ldots a^{\alpha_2}\ldots
b^{\beta_2}\ldots)$ when reading the word cyclically.
 \end {dfn}

 \begin {rk}
If we draw the chords of a chord diagram inside the circle then
linked chords intersect each other.
 \end {rk}

 \begin {dfn}
The {\em adjacency matrix} of a chord diagram $D$ with enumerated
$n$ chords is $n\times n$ matrix $A(D)=(a_{ij})$ over $\mathbb{Z}_2$
defined by
  \begin {enumerate}
   \item
$a_{ii}$ is the framing of the chord with the number $i$, i.e.\
either $G$ or $\pm1$;
   \item
$a_{ij}=1$, $i\ne j$ if and only if the chords with the numbers $i$
and $j$ are linked;
   \item
$a_{ij}=0$, $i\ne j$, if and only if $i$ and $j$ unlinked.
  \end {enumerate}
 \end {dfn}

 \begin {example}
Let $D$ be the framed chord diagram depicted in Fig.~\ref {frchdi}.
Enumerate all the chords of $D$: the chord $aa$ has the number $1$,
the chord $b$ has the number $2$ etc. Then
 $$
A(D)=\left(\begin {array}{ccccc}
 0 & 1 & 0 & 0 & 0\\
 1 & 0 & 1 & 0 & 1\\
 0 & 1 & 1 & 0 & 1\\
 0 & 0 & 0 & G & 1\\
 0 & 1 & 1 & 1 & 1
\end {array}\right).
 $$
 \end {example}

Assume we are given a chord diagram $D$ with all the chords having
framing $\pm1$ (with no Gaussian chord).

\begin {figure} \centering\includegraphics[width=400pt]{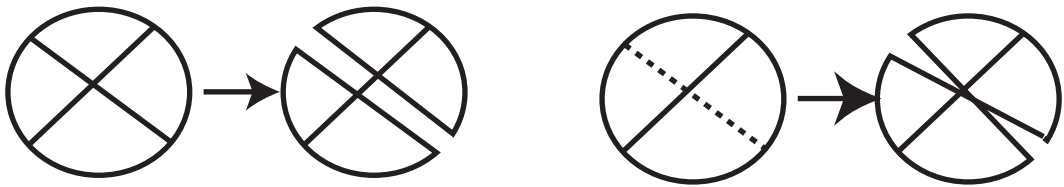}
  \caption{A Surgery of the Circuit along a Chord} \label{surger}
 \end {figure}

\begin {figure} \centering\includegraphics[width=200pt]{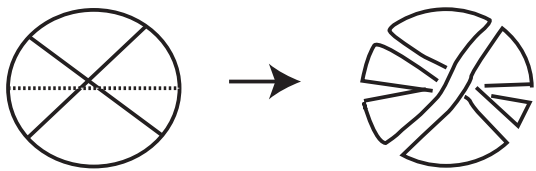}
  \caption{The Manifold $M(D)$} \label{result}
 \end {figure}

 \begin {dfn}
Define the {\em surgery along a set of chords} as follows. For every
chord having the framing $0$ (resp., $1$), we draw a parallel
(resp., intersecting) chord near it and remove the arc of the circle
between adjacent ends of the chords as in Fig.~\ref {surger}. By a
small perturbation, the picture in $\mathbb{R}^{2}$ is transformed
into a one-manifold in $\mathbb{R}^{3}$. This manifold $M(D)$ is the
{\em result of surgery}, see Fig.~\ref {result}.
 \end {dfn}

Surprisingly, the number of connected components of $M(D)$ can be
determined from the adjacency matrix $A(D)$ of $D$.

 \begin {theorem}[\cite {CL},~\cite {Sob},~\cite {BG},~\cite {Tr}]\label {th:sob}
Let $D$ be a chord diagram with all the chords having the framing
$\pm1$. Then the number of connected components of $M(D)$ equals
$\corank A(D)+1$, where $A(G)$ is the adjacency matrix of $D$ over
$\mathbb{Z}_2$ and $\corank$ is calculated over $\mathbb{Z}_2$.
  \end {theorem}

 \begin {rk}
All the matrices are considered over $\mathbb{Z}_2$ and we will not
indicate it explicitly.
 \end {rk}

By using Theorem~\ref {th:sob} we can formulate a criterion of the
existence of a Gauss circuit.

Let $D$ be a framed chord diagram with the adjacency matrix $A(D)$.
Construct the matrix $\widehat{A}(D)$ by deleting the rows and
columns of $A(D)$ corresponding to Gaussian chords.

 \begin {theorem}[\cite {IM1},~\cite {IM2}]\label{th:exga}
Let $H$ be a framed $4$-graph and $U$ be an Euler tour of $H$. Then
$H$ has a Gauss circuit if and only if
$\corank(\widehat{A}(D)+E)=0$, here $D$ is a framed chord diagram
constructed from $U$ and $E$ is the identity matrix.
 \end {theorem}

 \begin {proof}
The proof immediately follows from Fig.~\ref {smooth_ver}. In order
to get a Gauss circuit, i.e.\ the tour while traveling along it we
pass from $e_3$ to $e_1$, we have to delete all the Gaussian chords,
to replace all the chords having the framing $1$ with intersecting
chords and all the chords having the framing $0$ with parallel
chords.

 \begin {figure} \centering\includegraphics[width=250pt]{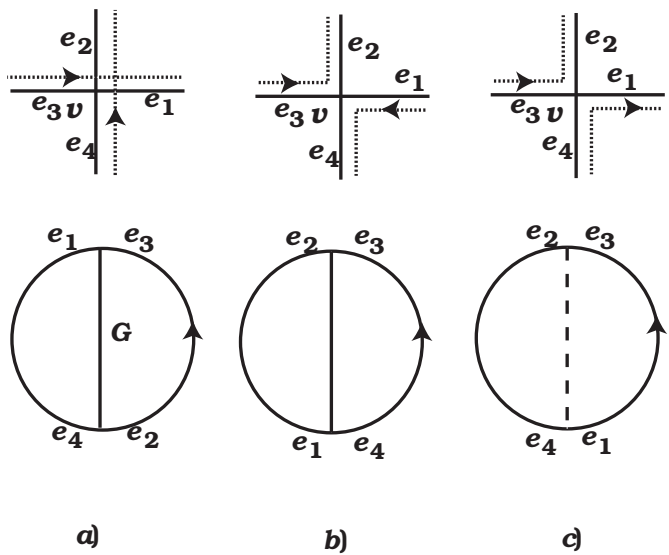}
  \caption{The Structure of Framed Chord Diagram} \label{smooth_ver}
 \end {figure}
 \end {proof}

\subsection {The Gauss Circuit}

Let $H$ be a framed $4$-graph having the Gauss circuit and let $U$
be an Euler tour of $H$. By using Corollary~\ref {cor:ro_ci} we can
assume that $m(U)$ (the corresponding chord diagram $D$) has no
Gaussian vertices ($U$ is a rotating circuit). The main result of
the whole paper is the following theorem.

 \begin {theorem}\label {th:ga_cir}
The adjacency matrix of the Gauss circuit is equal to
$(A(D)+E)^{-1}$ {\rom (}over $\mathbb{Z}_2${\rom)} up to diagonal
elements.
 \end {theorem}

 \begin {proof}
Let $V(H)=\{v_1,\dots,v_n\}$. It is not difficult to show that the
following two operations applied to $D$ decrease the number of
non-Gaussian chords:
  \begin {enumerate}
   \item
the framed star operation applying to a non-Gaussian chord having
the framing $0$;
   \item
$(((m*a)*b)*a)$, here $m$ is a framed double occurrence word,
$a,\,b$ are non-Gaussian letters (chords) having the framing $1$ and
they alternate in $m$ (are linked).
  \end {enumerate}
We call these operations {\em decreasing operations}. The decreasing
operations change an Euler tour $U$ on a $4$-graph, and the new
Euler tour has the number of non-Gaussian vertices smaller than $U$
has, see Fig.~\ref {croci} a),~\ref {AdGa}.

 \begin {figure} \centering\includegraphics[width=250pt]{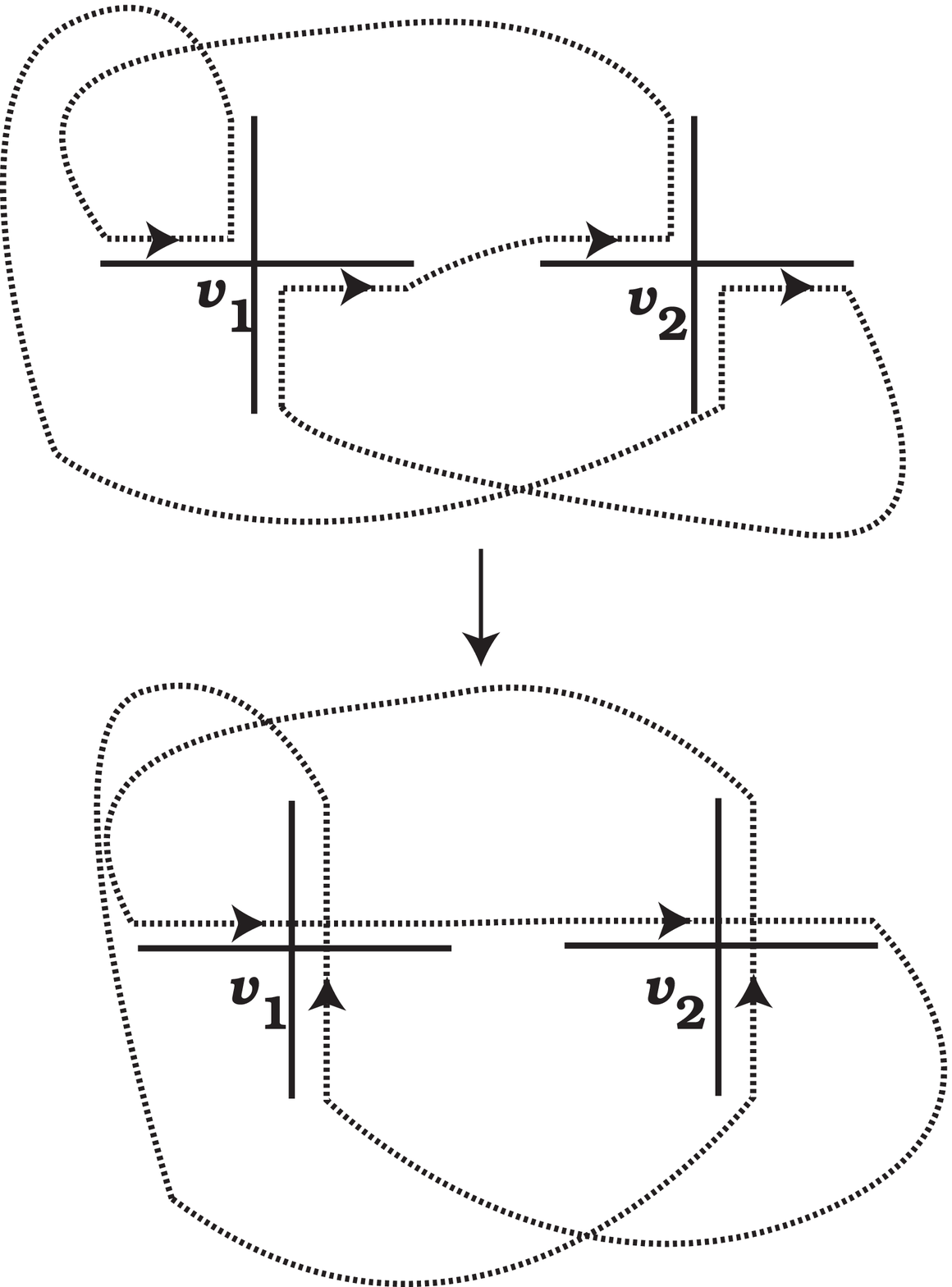}
  \caption{The Decreasing Operation} \label{AdGa}
 \end {figure}

Let $D$ be a framed chord diagram and let $A(D)$ be its adjacency
matrix. Let us apply decreasing operations. Without loss of
generality we may assume that the decreasing operations are applied
to the chords having the first numbers in our numeration. Then the
first decreasing operation is
 $$
A(D)=\left(\begin {array}{ccc} 0 & \0^\top & \1^\top \\
\0 & A_0 & A_1\\
\1 & A^\top_1 & A_2
 \end
 {array}\right)\rightsquigarrow
A(D')=\left(\begin {array}{ccc} G & \0^\top & \1^\top \\
\0 & A_0 & A_1\\
\1 & A^\top_1 & A_2+(1)
 \end{array}\right),
 $$
and the second one is
 $$
A(D)=\left(\begin {array}{cccccc}
1 & 1 & {\mathbf 0}^\top & {\mathbf 1}^\top & {\mathbf 0}^\top & {\mathbf 1}^\top\\
1 & 1 & {\mathbf 0}^\top & {\mathbf 0}^\top & {\mathbf 1}^\top & {\mathbf 1}^\top\\
{\mathbf 0} & {\mathbf 0} & A_0 & A_1 & A_2 & A_3\\
{\mathbf 1} & {\mathbf 0} & A^\top_1 & A_4 & A_5 & A_6\\
{\mathbf 0} & {\mathbf 1} & A^\top_2 & A^\top_5 & A_7 & A_8\\
{\mathbf 1} & {\mathbf 1} & A^\top_3 & A^\top_6 & A^\top_8 & A_9
 \end {array}\right)
 $$
 $$
\rightsquigarrow A(D')=\left(\begin {array}{cccccc}
G & 1 & {\mathbf 0}^\top & {\mathbf 0}^\top & {\mathbf 1}^\top & {\mathbf 1}^\top\\
1 & G & {\mathbf 0}^\top & {\mathbf 1}^\top & {\mathbf 0}^\top & {\mathbf 1}^\top\\
{\mathbf 0} & {\mathbf 0} & A_0 & A_1 & A_2 & A_3\\
{\mathbf 0} & {\mathbf 1} & A^\top_1 & A_4 & A_5+(1) & A_6+(1)\\
{\mathbf 1} & {\mathbf 0} & A^\top_2 & A^\top_5+(1) & A_7 & A_8+(1)\\
{\mathbf 1} & {\mathbf 1} & A^\top_3 & A^\top_6+(1) & A^\top_8+(1) &
A_9
 \end {array}\right),
 $$
where bold characters $\0$ and $\1$ indicate column vectors with all
entries the same $0$ and $1$, respectively, and $A_i$ are matrices.

We will successively apply these operations to $D$. On the next goal
is to show that after applying these two decreasing operations to
the framed chord diagram $D$ having no Gaussian vertices we will get
the framed chord diagram with the adjacency matrix $(A(D)+E)^{-1}$
up to diagonal elements.

To get the matrix $(A(D)+E)^{-1}$ we will perform elementary
manipulations with rows of $B(D)=A(D)+E$, $\det(A(D)+E)=1$. Let us
construct the matrix $(A(D)+E|E)$ with the size $n\times 2n$. We
denote by $\widehat{M}_{ij\ldots k}$ be the matrix obtained from $M$
by deleting $i,j,\ldots,k$-th rows and $i,j,\ldots,k$-th columns.

{\em The induction base}. As $\det B(D)=1$ then either there is a
diagonal element equal to $1$ or there are two diagonal elements
with the numbers $i$ and $j$ such that $b_{ii}=b_{ij}=0$,
$b_{ij}=b_{ji}=1$.

In the first case, without loss of generality assume that
$b_{11}=1$. Then after performing elementary manipulations with
$B(D)$ using the first row we get
 $$
B(D)=A(D)+E=\left(\begin {array}{ccc}
1 & \0^\top & \1^\top \\
\0 & A_0+E & A_1\\
\1 & A^\top_1 & A_2+E
 \end {array}\right)
 $$
 $$
\rightsquigarrow B'(D)=\left(\begin {array}{ccc} 1 & \0^\top & \1^\top \\
\0 & A_0+E & A_1\\
\0 & A^\top_1 & A_2+E+(1)
 \end
 {array}\right)
 $$
and
 $$
(B(D)|E)\rightsquigarrow(B'(D)|E')
 $$
 $$
=\left(\begin {array}{ccc} 1 & \0^\top & \1^\top \\
\0 & A_0+E & A_1\\
\0 & A^\top_1 & A_2+E+(1)
 \end {array}\right|\left.\begin {array}{ccc} 1 & \0^\top & \0^\top \\
\0 & E & 0\\
\1 & 0 & E \end {array}\right).
 $$
After performing the first decreasing operation to $D$ the chord
corresponding to $v_1$ becomes a Gaussian chord and the adjacencies
of non-Gaussian chords are defined by matrix $\widehat{B'}_1(D)$ and
the other adjacencies are defined by the first column of $E'$ (up to
diagonal elements).

In the second case, we may assume without loss of generality
$b_{11}=b_{22}=0$, $b_{12}=b_{21}=1$. Then after performing
elementary manipulations applied to the first two rows of $B(D)$, we
get
 $$
B(D)=A(D)+E=\left(\begin {array}{cccccc}
0 & 1 & {\mathbf 0}^\top & {\mathbf 1}^\top & {\mathbf 0}^\top & {\mathbf 1}^\top\\
1 & 0 & {\mathbf 0}^\top & {\mathbf 0}^\top & {\mathbf 1}^\top & {\mathbf 1}^\top\\
{\mathbf 0} & {\mathbf 0} & A_0+E & A_1 & A_2 & A_3\\
{\mathbf 1} & {\mathbf 0} & A^\top_1 & A_4+E & A_5 & A_6\\
{\mathbf 0} & {\mathbf 1} & A^\top_2 & A^\top_5 & A_7+E & A_8\\
{\mathbf 1} & {\mathbf 1} & A^\top_3 & A^\top_6 & A^\top_8 & A_9+E
 \end {array}\right)
 $$
 $$
\rightsquigarrow B'(D)=\left(\begin {array}{cccccc}
1 & 0 & {\mathbf 0}^\top & {\mathbf 0}^\top & {\mathbf 1}^\top & {\mathbf 1}^\top\\
0 & 1 & {\mathbf 0}^\top & {\mathbf 1}^\top & {\mathbf 0}^\top & {\mathbf 1}^\top\\
{\mathbf 0} & {\mathbf 0} & A_0+E & A_1 & A_2 & A_3\\
{\mathbf 0} & {\mathbf 0} & A^\top_1 & A_4+E & A_5+(1) & A_6+(1)\\
{\mathbf 0} & {\mathbf 0} & A^\top_2 & A^\top_5+(1) & A_7+E & A_8+(1)\\
{\mathbf 0} & {\mathbf 0} & A^\top_3 & A^\top_6+(1) & A^\top_8+(1) &
A_9+E
 \end {array}\right).
 $$
and
 $$
(B(D)|E)\rightsquigarrow(B'(D)|E')
 $$
 $$
=\left(\begin {array}{cccccc}
1 & 0 & {\mathbf 0}^\top & {\mathbf 0}^\top & {\mathbf 1}^\top & {\mathbf 1}^\top\\
0 & 1 & {\mathbf 0}^\top & {\mathbf 1}^\top & {\mathbf 0}^\top & {\mathbf 1}^\top\\
{\mathbf 0} & {\mathbf 0} & A_0+E & A_1 & A_2 & A_3\\
{\mathbf 0} & {\mathbf 0} & A^\top_1 & A_4+E & A_5+(1) & A_6+(1)\\
{\mathbf 0} & {\mathbf 0} & A^\top_2 & A^\top_5+(1) & A_7+E & A_8+(1)\\
{\mathbf 0} & {\mathbf 0} & A^\top_3 & A^\top_6+(1) & A^\top_8+(1) &
A_9+E
 \end {array}\right.
 $$
 $$
 \left|\begin {array}{cccccc}
0 & 1 & {\mathbf 0}^\top & {\mathbf 0}^\top & {\mathbf 0}^\top & {\mathbf 0}^\top\\
1 & 0 & {\mathbf 0}^\top & {\mathbf 0}^\top & {\mathbf 0}^\top & {\mathbf 0}^\top\\
{\mathbf 0} & {\mathbf 0} & E & 0 & 0 & 0\\
{\mathbf 0} & {\mathbf 1} & 0 & E & 0 & 0\\
{\mathbf 1} & {\mathbf 0} & 0 & 0 & E & 0\\
{\mathbf 1} & {\mathbf 1} & 0 & 0 & 0 & E
 \end {array}\right).
 $$
After performing the second decreasing operation to $D$ the chords
corresponding to $v_1$ and $v_2$ become Gaussian chords and the
adjacencies of the non-Gaussian chords are defined by matrix
$\widehat{B'}_{12}(D)$ and the other adjacencies are defined by the
first two columns of $E'$.

{\em The induction step}.  Let us suppose that we have performed $k$
decreasing operations. After these operations the matrix $(B(D)|E)$
is transformed into a matrix
 $$
(B'(D)|E')=\left(\begin {array}{cc} E & C\\ 0 & R\end
{array}\right|\left.\begin {array}{cc} F & 0\\ S & E\end
{array}\right)
 $$
and $F$ is a $l\times l$ matrix, $S$ is a symmetric matrix. Then the
new framed chord diagram contains $l$ Gaussian chords, and the
adjacencies of non-Gaussian chords are defined by $R$ and the other
adjacencies are defined by first $l$ rows of $E'$. As $\det B'(D)=1$
then $\det R=1$, and in the matrix $R$ there is either a diagonal
element equal to $1$ or there are numbers $p$ and $q$ such that
$r_{pp}=r_{qq}=0$, $r_{pq}=r_{qp}=1$.

Let us consider the first case. Without loss of generality we may
assume that $r_{11}=1$. In this case we apply the first decreasing
operation. We will get
 $$
(B'(D)|E')=\left(\begin {array}{cc} E & C\\ 0 & R\end
{array}\right|\left.\begin {array}{cc} F & 0\\ S & E\end
{array}\right)
 $$
 $$
=\left(\begin {array}{cccc} E & C_1 & C_2 & C_3\\
0 & 1 & \0^\top & \1^\top\\
\0 & \0 & R_1 & R_2\\
\0 & \1 & R^\top_2 & R_3
\end {array}\right|\left.\begin {array}{cccc} F & 0 & 0&0\\
S_1 & 1 & \0^\top & \0^\top\\
S_2 & \0 & E & 0\\
S_3 & \0 & 0 & E\end {array}\right)
 $$
 $$
\rightsquigarrow\left(\begin {array}{cccc} E & 0 & C'_2 & C'_3\\
0 & 1 & \0^\top & \1^\top\\
\0 & \0 & R_1 & R_2\\
\0 & \0 & R^\top_2 & R_3+(1)
\end {array}\right|\left.\begin {array}{cccc} F'_1 & F'_2 & 0& 0\\
S_1 & 1 & \0^\top & \0^\top\\
S_2 & \0 & E & 0\\
S'_3 & \1 & 0 & E\end {array}\right)
 $$
 $$
=\left(\begin {array}{cc} E & C'\\ 0 & R'\end
{array}\right|\left.\begin {array}{cc} F' & 0\\ S' & E\end
{array}\right)=(B''(D)|E''),
 $$
where $F'$ is a $(l+1)\times(l+1)$ matrix, $D'$ is a symmetric
matrix. The number of Gaussian vertices is $l+1$, and the
adjacencies of non-Gaussian vertices are defined by $R'$ and the
other adjacencies are defined by first $l$ rows of $E''$. The second
case is consider analogously to the first one.

We end up with the matrix
 $$
\left(\begin {array}{c} E\end {array}\right|\left.\begin {array}{c}
(A(D)+E)^{-1}\end {array}\right)
 $$
and the framed chord diagram having only Gaussian vertices. The
adjacency matrix of this chord diagram is $(A(D)+E)^{-1}$ up to
diagonal elements. We have proved Theorem for non-diagonal vertices.
But we know that all the diagonal elements is $G$.
 \end {proof}

 \begin {rk}
Let $U_1$ and $U_2$ be two rotating circuits, and let $D_1$ and
$D_2$ be its framed chord diagram such that $\det(A(D_i)+E)=1$. Then
the matrices $(A(D_1)+E)^{-1}$ and $(A(D_2)+E)^{-1}$ coincide only
up to diagonal elements.
 \end {rk}

\begin {figure} \centering\includegraphics[width=150pt]{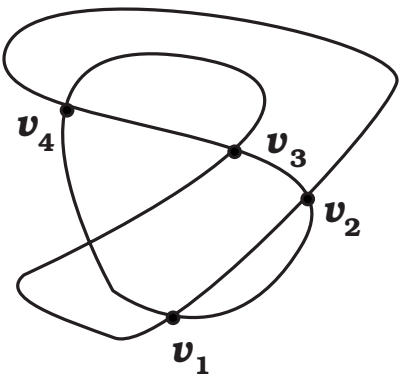}
  \caption{The Decreasing Operation} \label{examp}
 \end {figure}

 \begin {example}
Consider the framed $4$-graph having $4$ vertices $v_i$, Fig.~\ref
{examp}. Let $U_1$ and $U_2$ be two rotating circuits given by the
framed double occurrence cyclic words
$m(U_1)=(v_1v_4v_2v_1^{-1}v_2v_3v_4v_3)$ and
$m(U_2)=(v_1v_4v_3v_4v_2v_3v_1v_2^{-1})$, respectively. Then
 $$
A(m(U_1))=\left(\begin {array}{ccccc}
 1 & 1 & 0 & 1\\
 1 & 0 & 0 & 0\\
 0 & 0 & 0 & 1\\
 1 & 0 & 1 & 0
 \end {array}\right),\quad
A(m(U_2))=\left(\begin {array}{ccccc}
 0 & 1 & 0 & 0\\
 1 & 1 & 1 & 0\\
 0 & 1 & 0 & 1\\
 0 & 0 & 1 & 0
 \end {array}\right).
 $$
We get
 $$
(A(m(U_1))+E)^{-1}=\left(\begin {array}{ccccc}
 0 & 0 & 1 & 1\\
 0 & 1 & 1 & 1\\
 1 & 1 & 0 & 1\\
 1 & 1 & 1 & 1
 \end {array}\right),
 $$
 $$
(A(m(U_2))+E)^{-1}=\left(\begin {array}{ccccc}
 1 & 0 & 1 & 1\\
 0 & 0 & 1 & 1\\
 1 & 1 & 1 & 1\\
 1 & 1 & 1 & 0
 \end {array}\right)
 $$
and
 $$
A=\left(\begin {array}{ccccc}
 G & 0 & 1 & 1\\
 0 & G & 1 & 1\\
 1 & 1 & G & 1\\
 1 & 1 & 1 & G
 \end {array}\right)
 $$
is the adjacency matrix of Gauss circuit given by
$(v_1v_4v_3v_1v_2v_4v_3v_2)$.
 \end {example}

 \begin {dfn}
A chord diagram is called a {\em $d$-diagram} if its set of chords
can be split into two sets such that the chords from one set do not
intersect each other and all the chords are non-Gaussian having
framing $0$.
 \end {dfn}

The next Corollary immediately follows from the criterion of the
planarity of an immersed curve and from atom theory, see~\cite
{Fom,Ma3}.

 \begin {cor}[V.O.Manturov]\label {cor:man}
If $D$ is a chord diagram, all the chords of $D$ are non-Gaussian
having framing $0$ and $\det(A+E)=1$. Then for any $n$ numbers
$\lambda_1,\ldots,\lambda_n\in\mathbb{Z}_2$ such that
$\det\bigl((A+E)^{-1}+{\operatorname{diag}}(\lambda_1,\ldots,\lambda_n)\bigr)=1$
the matrix
$\bigl((A+E)^{-1}+{\operatorname{diag}}(\lambda_1,\ldots,\lambda_n)\bigr)^{-1}$
has $1$ on the diagonal. Moreover, if $D$ is a $d$-diagram then the
matrix
$\bigl((A+E)^{-1}+{\operatorname{diag}}(\lambda_1,\ldots,\lambda_n)\bigr)^{-1}$
is the adjacency matrix of a $d$-diagram. Here
${\operatorname{diag}}(\lambda_1,\ldots,\lambda_n)$ is the diagonal
matrix with $\lambda_1,\ldots,\lambda_n$ on its diagonal.
 \end {cor}


 \section {Adjacency Matrices}


It is well known that there are symmetric matrices over
$\mathbb{Z}_2$ which cannot be realised by chord diagrams, see~\cite
{Bou}, and adjacency matrices which can be realised by different
chord diagrams. In Theorem~\ref {th:ga_cir}, we have just used
elementary manipulations and adjacency matrices. It turns out that
Theorem~\ref {th:ga_cir} and Corollary~\ref {cor:man} can be
reformulated for arbitrary symmetric matrices.

In this section all the matrices are over $\mathbb{Z}_2$ with the
diagonal elements equal to $\pm1$. Consider the following operation
over the set of symmetric matrices. Let $A=(a_{ij})$ be a symmetric
$n\times n$ matrix. Let us fix an arbitrary element $a_{kk}=1$ and
construct the matrix $\Loc(A,k)=(\widetilde{a}_{ij})$, here
$\widetilde{a}_{pq}=a_{pq}+1$, $p,\,q\ne k$, if both $a_{pk}=1$ and
$a_{kq}=1$, and $\widetilde{a}_{pq}=a_{pq}$ otherwise. We call the
transformation $A\mapsto\Loc(A,k)$ a {\em local complementation of
the matrix $A$ at the element $a_{kk}$} (this operation is analogous
to the framed star).

It is not difficult to show that the matrices
$\Loc(\Loc(\Loc(A,i),j,i)$ and $\Loc(\Loc(\Loc(A,j),i,j)$ obtained
from a matrix $A$ coincide up to the diagonal elements with the
numbers $i,\,j$.

 \begin {dfn}
Let $A$ be a symmetric matrix with $a_{ii}=a_{jj}=0$,
$a_{ij}=a_{ji}=1$. A {\em pivot} operation is the transformation
$A\mapsto\widetilde{A}$, where the diagonal elements of
$\widetilde{A}$ coincide with the diagonal elements of $A$ and the
other elements of $\widetilde{A}$ coincide with the corresponding
elements of $\Loc(\Loc(\Loc(A,i),j,i)$.
 \end {dfn}

Let $\Sym(n,\mathbb{Z}_2)$ be the set of all symmetric $n\times n$
matrices over $\mathbb{Z}_2$. Consider two equivalence relations on
$\Sym(n,\mathbb{Z}_2)$. The first relation is: two matrices $A$ and
$B$ are said to be {\em equivalent up to diagonal}, denote this
equivalence relation by $A\sim_D B$, if $A$ and $B$ coincide up to
diagonal elements. The second equivalence is defined as follows: two
matrices $A$ and $B$ are said to be {\em obtained from each other by
changing circuit}, denote the second relation by $A\sim_C B$, if $A$
and $B$ are related by a sequence of local complementations and
pivot operations.

 \begin {lem}[\cite {IM1}, \cite {IM2}]
If $\det(A+E)=1$ and $B\sim_CA$ then $\det(B+E)=1$, here $E$ is the
identity matrix.
 \end {lem}

Let $\Sym_+(n,\mathbb{Z}_2)\subset\Sym(n,\mathbb{Z}_2)$ be the
subset of the set of symmetric matrices consisting of matrices $A$
with $\det(A+E)=1$.

 \begin {cor}
The relation $\sim_C$ is also the equivalence relation on
$\Sym_+(n,\mathbb{Z}_2)$.
 \end {cor}

Consider the two set
 $$
\mathfrak{L}(n)=\Sym(n,\mathbb{Z}_2)/\sim_D\qquad\hbox{and}\qquad
\mathfrak{G}(n)=\Sym_+(n,\mathbb{Z}_2)/\sim_{C}.
 $$

 \begin {lem}\label {lem:nondeg}
Every element of $\mathfrak{L}(n)$ has a representative with the
determinant equal to $1$.
 \end {lem}

 \begin {proof}
Let us prove this lemma by induction on the size of a matrix.

{\em The induction base}. For $n=1$ the claim of the lemma is
evident.

{\em The induction step}. Assume the statement of the lemma holds
for $n-1$ and let $A$ be a $n\times n$ matrix. By the induction
hypothesis, we can assume $\det A^{11}=1$, $A^{ij}$ is the algebraic
complement to $a_{ij}$. Then either
 $$
\det
A=a_{11}A^{11}+\sum\limits_{j=2}^na_{1j}A^{1j}=a_{11}+\sum\limits_{j=2}^na_{1j}A^{1j}=1
 $$
or
 $$
\det\widetilde{A}=(a_{11}+1)A^{11}+\sum\limits_{j=2}^na_{1j}A^{1j}=a_{11}+1+\sum\limits_{j=2}^na_{1j}A^{1j}=1,
 $$
where the matrix $\widetilde{A}$ is different from $A$ only by the
element $\widetilde{a}_{11}$.
 \end {proof}

 \begin {lem}\label {lem:eq1}
Let $B$ and $\widetilde{B}$ be two matrices over $\mathbb{Z}_2$ with
$\det B=\det\widetilde{B}=1$, and let $B$ and $\widetilde{B}$
coincide up to one element on the diagonal. Then the matrices
$B^{-1}+E$ and $\widetilde{B}^{-1}+E$ are related by a local
complementation.
 \end {lem}

 \begin {proof}
Without loss of generality we may assume that $B=(b_{ij})$,
$\widetilde{B}=(\widetilde{b}_{ij})$ are $n\times n$ matrices and
$\widetilde{b}_{nn}=b_{nn}+1=0$, $\widetilde{b}_{ij}=b_{ij}$, $i\ne
n$ or $j\ne n$. We will perform elementary manipulations with rows
of $B$ and $\widetilde{B}$ to get two identity matrices. Then we
will apply these elementary manipulations to two identity matrix to
get the inverse matrices.

Using the equality $\det B=\det\widetilde{B}=1$ we have
 $$
1=\det\widetilde{B}=\det B+\det \widehat{B}^n_n=1+\det
\widehat{B}^n_n,\quad \det\widehat{B}^n_n=0,
 $$
 $$
\rank B=n,\quad \rank \widehat{B}^n_n=n-2,
 $$
here $\widehat{B}^n_n$ is the matrix obtained from $B$ by deleting
the $n$-th row and $n$-th column. As $\widehat{B}^n_n$ is a
symmetric matrix, without loss of generality we may assume that
$\det C=1$, here $C$ is the matrix obtained from $\widehat{B}^n_n$
by deleting the $(n-1)$-th row and $(n-1)$-th column.

Performing elementary manipulations with rows of $B$ and
$\widetilde{B}$ (the first $(n-1)$ rows of $\widetilde{B}$ are the
same as the ones of $B$), we get
 $$
B\rightsquigarrow\left(\begin {array}{ccc}
E       & \u & \v\\
\0^\top &  0 &  1\\
\0^\top &  1 &  1 \end {array}\right|\left.\begin {array}{ccc}
F  & \0 & \0\\
\a^\top &  1 &  0\\
\b^\top &  0 &  1
 \end {array}\right),
 $$
 $$
\widetilde{B}\rightsquigarrow\left(\begin {array}{ccc}
E       & \u & \v\\
\0^\top &  0 &  1\\
\0^\top &  1 &  0 \end {array}\right|\left.\begin {array}{ccc}
F  & \0 & \0\\
\a^\top &  1 &  0\\
\b^\top &  0 &  1
 \end {array}\right),
 $$
here $F$ is a $(n-2)\times(n-2)$-matrix, $\a,\,\b,\,\u,\,\v$ are
$(n-2)$-column vectors. Further, performing elementary manipulations
with rows, we have
 $$
\widetilde{B}\rightsquigarrow\left(\begin {array}{ccc}
E       & \u & \v\\
\0^\top &  0 &  1\\
\0^\top &  1 &  0 \end {array}\right|\left.\begin {array}{ccc}
F       & \0 & \0\\
\a^\top &  1 &  0\\
\b^\top &  0 &  1
 \end {array}\right)\rightsquigarrow\left(\begin {array}{ccc}
E       & \u & \v\\
\0^\top &  1 &  0\\
\0^\top &  0 &  1 \end {array}\right|\left.\begin {array}{ccc}
F  & \0 & \0\\
\b^\top &  0 &  1\\
\a^\top &  1 &  0
 \end {array}\right)
 $$
 $$
\rightsquigarrow\left(\begin {array}{ccc}
E       & \u & \0\\
\0^\top &  1 &  0\\
\0^\top &  0 &  1 \end {array}\right|\left.\begin {array}{ccc}
F_1     & \b & \0\\
\b^\top &  0 &  1\\
\a^\top &  1 &  0
 \end {array}\right)
\rightsquigarrow\left(\begin {array}{ccc}
E       & \0 & \0\\
\0^\top &  1 &  0\\
\0^\top &  0 &  1 \end {array}\right|\left.\begin {array}{ccc}
F_2     & \b & \a\\
\b^\top &  0 &  1\\
\a^\top &  1 &  0
 \end {array}\right),
 $$
i.e.\ $\u=\a,\,\v=\b$ (the inverse matrix to a symmetric matrix is
symmetric),
 $$
B\rightsquigarrow\left(\begin {array}{ccc}
E       & \u & \v\\
\0^\top &  0 &  1\\
\0^\top &  1 &  1 \end {array}\right|\left.\begin {array}{ccc}
F  & \0 & \0\\
\a^\top &  1 &  0\\
\b^\top &  0 &  1
 \end {array}\right)
 $$
 $$
\rightsquigarrow
 \left(\begin {array}{ccc}
E       & \a & \b\\
\0^\top &  1 &  0\\
\0^\top &  0 &  1 \end {array}\right|\left.\begin {array}{ccc}
F  & \0 & \0\\
\a^\top+\b^\top &  1 &  1\\
\a^\top &  1 &  0
 \end {array}\right)
 $$
 $$
\rightsquigarrow
 \left(\begin {array}{ccc}
E       & \a & \0\\
\0^\top &  1 &  0\\
\0^\top &  0 &  1 \end {array}\right|\left.\begin {array}{ccc}
F_1 & \b & \0\\
\a^\top+\b^\top &  1 &  1\\
\a^\top &  1 &  0
 \end {array}\right)
 $$
 $$
 \rightsquigarrow
 \left(\begin {array}{ccc}
E       & \0 & \0\\
\0^\top &  1 &  0\\
\0^\top &  0 &  1 \end {array}\right|\left.\begin {array}{ccc}
F_3 & \a+\b & \a\\
\a^\top+\b^\top &  1 &  1\\
\a^\top &  1 &  0
 \end {array}\right).
 $$
It is not difficult to see that $F_3$ is obtained from $F_2$ by
adding $\a^\top$ to the rows of $F_2$ corresponding to the rows of
$B^{-1}$ which have the last element equal to $1$. Therefore, the
matrix $B^{-1}+E$ is obtained from $\widetilde{B}^{-1}+E$ by the
local complementation at the element corresponding to
$\widetilde{b}_{nn}$.
 \end {proof}

 \begin {lem}\label {lem:eq2}
Let $B$ and $\widetilde{B}$ be two matrices over $\mathbb{Z}_2$ with
$\det B=\det\widetilde{B}=1$, and let $B$ and $\widetilde{B}$
coincide up to two elements on the diagonal with numbers $i$ and
$j$. Suppose that $\det\widehat{B}^i_i=\det\widehat{B}^j_j=1$, here
$\widehat{B}^k_k$ is the matrix obtained from $B$ by deleting the
$k$-th row and $k$-th column. Then the matrices $B^{-1}+E$ and
$\widetilde{B}^{-1}+E$ are related by a pivot operation.
 \end {lem}

 \begin {proof}
Without loss of generality, we may assume that $B=(b_{ij})$,
$\widetilde{B}=(\widetilde{b}_{ij})$ are $n\times n$ matrices and
$\widetilde{b}_{(n-1)(n-1)}=b_{(n-1)(n-1)}+1$,
$\widetilde{b}_{nn}=b_{nn}+1$, $\widetilde{b}_{ij}=b_{ij}$ for
$(i,j)\ne(n-1,n-1),(n,n)$. We will perform elementary manipulations
with rows of $B$ and $\widetilde{B}$ to get two identity matrices.
Then we will apply these manipulations to two identity matrices to
get the inverse matrices.

By using the equality $\det B=\det\widetilde{B}=1$, we have
 $$
1=\det\widetilde{B}=\det B+\det \widehat{B}^{n-1}_{n-1}+\det
\widehat{B}^n_n+\det\widehat{B}^{(n-1)n}_{(n-1)n}=1+\det
\widehat{B}^{(n-1)n}_{(n-1)n},
 $$
 $$
\det\widehat{B}^{(n-1)n}_{(n-1)n}=0,\quad \rank
\widehat{B}^{n-1}_{n-1}=\rank \widehat{B}^n_n=n-1,\quad \rank
\widehat{B}^{(n-1)n}_{(n-1)n}=n-3,
 $$
here $\widehat{B}^{(n-1)n}_{(n-1)n}$ is the matrix obtained from
$\widehat{B}^n_n$ by deleting the $(n-1)$-th row and $(n-1)$-th
column. As $\widehat{B}^{(n-1)n}_{(n-1)n}$ is a symmetric matrix,
without loss of generality we may assume that $\det C=1$, here $C$
is the matrix obtained from $\widehat{B}^{(n-1)n}_{(n-1)n}$ by
deleting the $(n-2)$-th row and $(n-2)$-th column. It is not
difficult to show that the matrices obtained from $\widetilde{B}$ by
deleting the $n$-th row, $n$-th column and the $(n-1)$-th row, the
$(n-1)$-th column, respectively, are both nondegenerate.

Performing elementary manipulations with rows of $B$ and
$\widetilde{B}$ (the first $(n-2)$ rows of $\widetilde{B}$ are the
same as the ones of $B$), we get
 $$
B\rightsquigarrow\left(\begin {array}{cccc}
E       & \u & \v & \w\\
\0^\top &  0 &  1 & 1\\
\0^\top &  1 &  1 & l\\
\0^\top &  1 &  l & 0\end {array}\right|\left.\begin {array}{cccc}
F       & \0 & \0 & \0\\
\a^\top &  1 &  0 & 0\\
\b^\top &  0 &  1 & 0\\
\c^\top &  0 &  0 & 1
 \end {array}\right),
 $$
 $$
\widetilde{B}\rightsquigarrow\left(\begin {array}{cccc}
E       & \u & \v & \w\\
\0^\top &  0 &  1 & 1\\
\0^\top &  1 &  0 & l\\
\0^\top &  1 &  l & 1\end {array}\right|\left.\begin {array}{cccc}
F       & \0 & \0 & \0\\
\a^\top &  1 &  0 & 0\\
\b^\top &  0 &  1 & 0\\
\c^\top &  0 &  0 & 1
 \end {array}\right),
 $$
here $F$ is a $(n-3)\times(n-3)$-matrix,
$\a,\,\b,\,\c,\,\u,\,\v,\,\w$ are $(n-3)$-column vectors, and
$l\in\{0,1\}$. Further, performing elementary manipulations with
rows, we have
 $$
B\rightsquigarrow\left(\begin {array}{cccc}
E       & \u & \v & \w\\
\0^\top &  0 &  1 & 1\\
\0^\top &  1 &  1 & l\\
\0^\top &  1 &  l & 0\end {array}\right|\left.\begin {array}{cccc}
F       & \0 & \0 & \0\\
\a^\top &  1 &  0 & 0\\
\b^\top &  0 &  1 & 0\\
\c^\top &  0 &  0 & 1
 \end {array}\right)
  $$
  $$
\rightsquigarrow \left(\begin {array}{cccc}
E       & \u & \v & \w\\
\0^\top &  1 &  0 & 0\\
\0^\top &  0 &  1 & 0\\
\0^\top &  0 &  0 & 1\end {array}\right|\left.\begin {array}{cccc}
F       & \0 & \0 & \0\\
l\a^\top+l\b^\top+(1+l)\c^\top &  l &  l & 1+l\\
l\a^\top+\b^\top+\c^\top &  l &  1 & 1\\
(1+l)\a^\top+\b^\top+\c^\top &  1+l &  1 & 1
 \end {array}\right)
 $$
 $$
\rightsquigarrow \Bigl(\begin {array}{c} E \end {array}\Big|
\begin {array}{c} B^{-1} \end {array}\Bigr),
 $$
here
  $$
B^{-1}=
 $$
 $$
\left(\begin {array}{cccc}
F_1       & l(\a+\b+\c)+\c & l\a+\b+\c & (1+l)\a+\b+\c\\
l\a^\top+l\b^\top+(1+l)\c^\top &  l &  l & 1+l\\
l\a^\top+\b^\top+\c^\top &  l &  1 & 1\\
(1+l)\a^\top+\b^\top+\c^\top &  1+l &  1 & 1
 \end {array}\right),
 $$
and
 $$
\widetilde{B}\rightsquigarrow\left(\begin {array}{cccc}
E       & \u & \v & \w\\
\0^\top &  0 &  1 & 1\\
\0^\top &  1 &  0 & l\\
\0^\top &  1 &  l & 1\end {array}\right|\left.\begin {array}{cccc}
F       & \0 & \0 & \0\\
\a^\top &  1 &  0 & 0\\
\b^\top &  0 &  1 & 0\\
\c^\top &  0 &  0 & 1
 \end {array}\right)
 $$
 $$
\rightsquigarrow \left(\begin {array}{cccc}
E       & \u & \v & \w\\
\0^\top &  1 &  0 & 0\\
\0^\top &  0 &  1 & 0\\
\0^\top &  0 &  0 & 1\end {array}\right|\left.\begin {array}{cccc}
F       & \0 & \0 & \0\\
l(\a^\top+\b^\top+\c^\top)+\b^\top &  l &  1+l & l\\
(1+l)\a^\top+\b^\top+\c^\top &  1+l &  1 & 1\\
l\a^\top+\b^\top+\c^\top & l &  1 & 1
 \end {array}\right)
 $$
 $$
\rightsquigarrow \Bigl(\begin {array}{c} E \end {array}\Big|
\begin {array}{c} \widetilde{B}^{-1} \end {array}\Bigr),
 $$
here
 $$
\widetilde{B}^{-1}=
 $$
 $$
\left(\begin {array}{cccc}
{F}_2       & l(\a+\b+\c)+\b & (1+l)\a+\b+\c & l\a+\b+\c\\
l(\a^\top+\b^\top+\c^\top)+\b^\top &  l &  1+l & l\\
(1+l)\a^\top+\b^\top+\c^\top &  1+l &  1 & 1\\
l\a^\top+\b^\top+\c^\top & l &  1 & 1
 \end {array}\right).
 $$

Let us investigate the matrices $F_1$ and $F_2$. We have $4$ cases.

a) If a row of the matrix $B^{-1}$ has the last two elements equal
to $0$ then the corresponding row of $\widetilde{B}^{-1}$ has also
the last two elements equal to $0$. We have two options: either the
rows of $F_1$ and $F_2$ are obtained from $F$ by adding
$l\a^\top+\b^\top+\c^\top$ and $(1+l)\a^\top+\b^\top+\c^\top$ to the
corresponding row of $F$ or they are equal to the corresponding row
of $F$. In both cases, we have the equality of rows of $B^{-1}$ and
$\widetilde{B}^{-1}$ having the last two elements equal to $0$.

b) If a row of the matrix $B^{-1}$ has the last two elements equal
to $1$ then the corresponding row of $\widetilde{B}^{-1}$ has also
the last two elements equal to $1$. We have two options: either the
rows of $F_1$ and $F_2$ are obtained from $F$ by adding
$l\a^\top+\b^\top+\c^\top$ or by adding
$(1+l)\a^\top+\b^\top+\c^\top$ to the corresponding row of $F$. In
both cases, we have the equality of rows of $B^{-1}$ and
$\widetilde{B}^{-1}$ having the last two elements equal to $1$.

c) If a row of the matrix $B^{-1}$ has the penultimate element equal
to $0$ and the last one is $1$ then the corresponding row of
$\widetilde{B}^{-1}$ has the penultimate element equal to $1$ and
the last one is $0$. We have two options: either the rows of $F_1$
and $F_2$ are obtained from $F$ by adding
$l\a^\top+l\b^\top+(1+l)\c^\top$ and $l(l\a^\top+\b^\top+\c^\top)$
for $F_1$ and $l(\a^\top+\b^\top+\c^\top)+\b^\top$ and
$l((1+l)\a^\top+\b^\top+\c^\top)$ for $F_2$ or by adding
$l\a^\top+l\b^\top+(1+l)\c^\top$, $(1+l)(l\a^\top+\b^\top+\c^\top)$
and $(1+l)\a^\top+\b^\top+\c^\top$ for $F_1$ and
$l(\a^\top+\b^\top+\c^\top)+\b^\top$,
$(1+l)((1+l)\a^\top+\b^\top+\c^\top)$ and $l\a^\top+\b^\top+\c^\top$
for $F_2$ to the corresponding row of $F$. In both cases, the sum of
the rows of $B^{-1}$ and $\widetilde{B}^{-1}$ is
$l\a^\top+\b^\top+\c^\top$.

d) If a row of the matrix $B^{-1}$ has has the penultimate element
equal to $1$ and the last one is $0$ then the corresponding row of
$\widetilde{B}^{-1}$ has the penultimate element equal to $0$ and
the last one is $1$. We have two options: either the rows of $F_1$
and $F_2$ are obtained from $F$ by adding
$l\a^\top+l\b^\top+(1+l)\c^\top$ and
$(1+l)(l\a^\top+\b^\top+\c^\top)$ for $F_1$ and
$l(\a^\top+\b^\top+\c^\top)+\b^\top$ and
$(1+l)((1+l)\a^\top+\b^\top+\c^\top)$ for $F_2$ or by adding
$l\a^\top+l\b^\top+(1+l)\c^\top$, $l(l\a^\top+\b^\top+\c^\top)$ and
$(1+l)\a^\top+\b^\top+\c^\top$ for $F_1$ and
$l(\a^\top+\b^\top+\c^\top)+\b^\top$,
$l((1+l)\a^\top+\b^\top+\c^\top)$ and $l\a^\top+\b^\top+\c^\top$ for
$F_2$ to the corresponding row of $F$. In both cases, the sum of the
rows of $B^{-1}$ and $\widetilde{B}^{-1}$ is
$(1+l)\a^\top+\b^\top+\c^\top$.

Therefore, the matrices $B^{-1}+E$ and $\widetilde{B}^{-1}+E$ are
related by a pivot operation.
 \end {proof}

 \begin {theorem}
  \begin {enumerate}
   \item
The map $\chi\colon\mathfrak{G}(n)\to\mathfrak{L}(n)$ given by the
formula $\chi[A]_C=[(A+E)^{-1}]_D$ is well defined.
   \item
There exists the inverse map
$\chi^{-1}\colon\mathfrak{L}(n)\to\mathfrak{G}(n)$.
  \end {enumerate}
 \end {theorem}

 \begin {proof}
Let $E_{ij}$ be the matrix with $1$ on the whole diagonal and the
element in the intersection of the $i$-th row and $j$-th column is
$1$, the others are $0$.

1) Let $A\sim_C\widetilde{A}$.

If $A$ and $\widetilde{A}$ are related by a pivot operation for the
first two elements. Then
 $$
B=A+E=\left(\begin {array}{cccccc}
1 & 1 & {\mathbf 0}^\top & {\mathbf 1}^\top & {\mathbf 0}^\top & {\mathbf 1}^\top\\
1 & 1 & {\mathbf 0}^\top & {\mathbf 0}^\top & {\mathbf 1}^\top & {\mathbf 1}^\top\\
{\mathbf 0} & {\mathbf 0} & A_0+E & A_1 & A_2 & A_3\\
{\mathbf 1} & {\mathbf 0} & A^\top_1 & A_4+E & A_5 & A_6\\
{\mathbf 0} & {\mathbf 1} & A^\top_2 & A^\top_5 & A_7+E & A_8\\
{\mathbf 1} & {\mathbf 1} & A^\top_3 & A^\top_6 & A^\top_8 & A_9+E
 \end {array}\right)
 $$
 $$
\widetilde{B}=\widetilde{A}+E=\left(\begin {array}{cccccc}
1 & 1 & {\mathbf 0}^\top & {\mathbf 0}^\top & {\mathbf 1}^\top & {\mathbf 1}^\top\\
1 & 1 & {\mathbf 0}^\top & {\mathbf 1}^\top & {\mathbf 0}^\top & {\mathbf 1}^\top\\
{\mathbf 0} & {\mathbf 0} & A_0+E & A_1 & A_2 & A_3\\
{\mathbf 0} & {\mathbf 1} & A^\top_1 & A_4+E & A_5+(1) & A_6+(1)\\
{\mathbf 1} & {\mathbf 0} & A^\top_2 & A^\top_5+(1) & A_7+E & A_8+(1)\\
{\mathbf 1} & {\mathbf 1} & A^\top_3 & A^\top_6+(1) & A^\top_8+(1) &
A_9+E
 \end {array}\right)
 $$
 $$
=BE_{1k_1}\ldots E_{1k_p}E_{2(k_p+1)}\ldots
E_{2k_q}E_{1(k_q+1)}\ldots E_{1n}
 $$
 $$
\cdot E_{2(k_q+1)}\ldots E_{2n}E_{12}E_{21}E_{12}=BM,
 $$
here $k_1>2,\ldots,k_p$ are the numbers of those columns which have
$1$ in the first row and $0$ in the second row, $k_p+1,\ldots,k_q$
are the numbers of those columns which have $0$ in the first row and
$1$ in the second row, and $k_q+1,\ldots,n$ are the numbers of those
columns which have $1$ in the first two rows.

We get $\widetilde{B}^{-1}=M^{-1}B^{-1}$. The last matrix is
obtained from $B^{-1}$ by adding rows to the first and second rows
of it. As matrices $\widetilde{B}^{-1}$ and $B^{-1}$ are symmetric
then $\widetilde{B}^{-1}$ might differ from $B^{-1}$ only by the
four elements located in the first two rows and columns. So we have
to prove the equality of $b^{12}=\widetilde{b}^{12}$,
$B^{-1}=(b^{ij})$, $\widetilde{B}^{-1}=(\widetilde{b}^{ij})$. We
have
 $$
b^{12}=\det\left(\begin {array}{ccccc}
1 & {\mathbf 0}^\top & {\mathbf 0}^\top & {\mathbf 1}^\top & {\mathbf 1}^\top\\
{\mathbf 0} & A_0+E & A_1 & A_2 & A_3\\
{\mathbf 1} & A^\top_1 & A_4+E & A_5 & A_6\\
{\mathbf 0} & A^\top_2 & A^\top_5 & A_7+E & A_8\\
{\mathbf 1} & A^\top_3 & A^\top_6 & A^\top_8 & A_9+E
 \end {array}\right)
 $$
 $$
=\det\left(\begin {array}{ccccc}
1 & {\mathbf 0}^\top & {\mathbf 0}^\top & {\mathbf 1}^\top & {\mathbf 1}^\top\\
{\mathbf 0} & A_0+E & A_1 & A_2 & A_3\\
{\mathbf 0} & A^\top_1 & A_4+E & A_5+(1) & A_6+(1)\\
{\mathbf 0} & A^\top_2 & A^\top_5 & A_7+E & A_8\\
{\mathbf 0} & A^\top_3 & A^\top_6 & A^\top_8+(1) & A_9+E+(1)
 \end {array}\right)
 $$
 $$
=\det\left(\begin {array}{cccc}
A_0+E & A_1 & A_2 & A_3\\
A^\top_1 & A_4+E & A_5+(1) & A_6+(1)\\
A^\top_2 & A^\top_5 & A_7+E & A_8\\
A^\top_3 & A^\top_6 & A^\top_8+(1) & A_9+E+(1)
 \end {array}\right),
 $$
 $$
\widetilde{b}^{12}=\det\left(\begin {array}{cccccc}
1 & {\mathbf 0}^\top & {\mathbf 1}^\top & {\mathbf 0}^\top & {\mathbf 1}^\top\\
{\mathbf 0} & A_0+E & A_1 & A_2 & A_3\\
{\mathbf 0} & A^\top_1 & A_4+E & A_5+(1) & A_6+(1)\\
{\mathbf 1} & A^\top_2 & A^\top_5+(1) & A_7+E & A_8+(1)\\
{\mathbf 1} & A^\top_3 & A^\top_6+(1) & A^\top_8+(1) & A_9+E
 \end {array}\right)
 $$
 $$
=\det\left(\begin {array}{cccccc}
1 & {\mathbf 0}^\top & {\mathbf 1}^\top & {\mathbf 0}^\top & {\mathbf 1}^\top\\
{\mathbf 0} & A_0+E & A_1 & A_2 & A_3\\
{\mathbf 0} & A^\top_1 & A_4+E & A_5+(1) & A_6+(1)\\
{\mathbf 0} & A^\top_2 & A^\top_5 & A_7+E & A_8\\
{\mathbf 0} & A^\top_3 & A^\top_6 & A^\top_8+(1) & A_9+E+(1)
 \end {array}\right)
 $$
 $$
=\det\left(\begin {array}{cccccc}
A_0+E & A_1 & A_2 & A_3\\
A^\top_1 & A_4+E & A_5+(1) & A_6+(1)\\
A^\top_2 & A^\top_5 & A_7+E & A_8\\
A^\top_3 & A^\top_6 & A^\top_8+(1) & A_9+E+(1)
 \end {array}\right)=b^{12}.
 $$
We have proven $B^{-1}\sim_D\widetilde{B}^{-1}$.

If $A$ and $\widetilde{A}$ are related by the local complementation
for the first element, then
 $$
B=A+E=\left(\begin {array}{ccc}
0 & \0^\top & \1^\top \\
\0 & A_0+E & A_1\\
\1 & A^\top_1 & A_2+E
\end {array}\right)
 $$
 $$
\widetilde{B}=\left(\begin {array}{ccc}
0 & \0^\top & \1^\top \\
\0 & A_0+E & A_1\\
\1 & A^\top_1 & A_2+(1)+E
\end {array}\right)=(A(G_1)+E)E_{1m}E_{1(m+1)}\ldots E_{1n},
 $$
here the numbers $m,m+1,\ldots,n$ correspond to the numbers of
columns containing $1$ in the first row.

We get $\widetilde{B}^{-1}=E_{1n}\ldots E_{1m} B^{-1}$. The matrix
$\widetilde{B}^{-1}$ is obtained from $B^{-1}$ by adding the rows
with numbers from $m$ to $n$ to the first row of it. As matrices
$\widetilde{B}^{-1}$ and $B^{-1}$ are symmetric then
$\widetilde{B}^{-1}$ might differ from $B^{-1}$ only by the first
diagonal element. So we have proven
$B^{-1}\sim_D\widetilde{B}^{-1}$.

If $A$ and $\widetilde{A}$ are related by pivot operations and local
complementations then consequently by applying two preceding cases
we get $(A+E)^{-1}\sim_D(\widetilde{A}+E)^{-1}$.

2) If $B\sim_D\widetilde{B}$ and $\det B=\det\widetilde{B}=1$ then,
by using Lemmas~\ref {lem:eq1} and~\ref {lem:eq2},
$B^{-1}+E\sim_C\widetilde{B}^{-1}+E$. By using Lemma~\ref
{lem:nondeg}, we see that there exists some $B$ with $\det B=1$ in
each class $[C]_D$. So we can define the inverse map
$\chi^{-1}\colon\mathfrak{L}(n)\to\mathfrak{G}(n)$ by
$\chi^{-1}([C]_D)=[B^{-1}+E]_C$.
 \end {proof}

 \begin {rk}
The map $\chi$ gives rise to an isomorphism between the set of
looped interlacement graphs modulo the Reidemeister moves, see~\cite
{TZ}, and graph-knots, see~\cite {IM1,IM2}. We shall address this
question in a sequel of the present paper.
 \end {rk}

 \section*{Acknowledgments}
The author is grateful to V.\,O.~Manturov, A.\,T.~Fomenko,
D.\,M.~Afanasiev, I.\,M.~Nikonov for their interest to this work and
to L.~Traldi for his discussion of~\cite {TZ}.

 \end {document}